\documentclass[11pt]{article}
\catcode`\@=11 \catcode`\@=12

\usepackage{amsmath}
\usepackage{amssymb}
\usepackage{amsthm}
\usepackage{oldgerm}

\usepackage[pdftex,colorlinks,urlcolor=blue,pdfstartview=FitH]{hyperref}

\newcommand{\Rm}{\mathbb{R}}

\newtheorem{lem}{Lemma}
\newtheorem{thm}{Theorem}

\def\proof {\noindent{\sc{Proof. }}}
\def\qed {\mbox{}\hfill {\small \fbox{}} \\}
\def\lto{\longrightarrow}
\def\lmto{\longmapsto}

\def\leq{\leqslant}
\def\geq{\geqslant}

\title{Lasry-Lions regularization and a Lemma of Ilmanen}
\author{Patrick  Bernard}
\date{ December 2009}
\begin{document}

\maketitle

Let $H$ be a Hilbert space. We define the following $\inf$ ($\sup$) convolution operators
acting on bounded functions $u:H\lto \Rm$:
$$
T_tu(x):=\inf_y \big( u(y)+\frac{1}{t}\|y-x\|^2\big)
$$
and 
$$
\check T_tu(x):= \sup_y \big( u(y)-\frac{1}{t}\|y-x\|^2\big).
$$
We have the relation
$$
T_t(-u)=-\check T_t(u).
$$
Recall that these operators form semi-groups, in the sense that 
$$
T_t\circ T_s=T_{t+s}\quad \text{and}\quad
\check T_t\circ \check T_s=\check T_{t+s}
$$
for all $t\geq 0$ and $s\geq 0$, as can be checked by direct calculation.
Note also that
$$\inf u\leq T_tu(x)\leq u(x) \leq \check T_tu(x)\leq \sup u
$$ 
for each $t\geq 0$ and each $x\in H$.
A   function $u:H\lto \Rm$ is called $k$-semi-concave, $k> 0$,
if the function $x\lto u(x)-\|x\|^2/k$ is concave. The function $u$
is called $k$-semi-convex if $-u$ is $k$-semi-concave.
A bounded function $u$ 
is $t$-semi-concave if and only if it belongs to the image of the operator $T_t$, this
follows from Lemma \ref{1} and Lemma \ref{3} below.
A function is called semi-concave if it is $k$-semi-concave for some $k>0$.
A function $u$ is said $C^{1,1}$ if it is Frechet differentiable and if the gradient
of $u$ is Lipschitz.
Note that a continuous function is $C^{1,1}$ if and only if it is semi-concave and semi-convex, 
see Lemma \ref{F}. 
Let us recall two important results in that language:

\begin{thm}(Lasry-Lions, \cite{LL})\label{LL}
Let $u$ be a bounded function.
For $0<s<t$, the function $\check T_s \circ T_t u$ is $C^{1,1}$ and,
if $u$ is uniformly continuous, then  it converges uniformly to
$u$ when $t\lto 0$.
\end{thm}

\begin{thm}\label{Il}(Ilmanen, \cite{Il})
Let $u\geq v$ be two bounded  functions on $H$ such that $u$ and $-v$ are semi-concave.
Then there exists a $C^{1,1}$ function $w$ such that $u\geq w\geq v$.
\end{thm}

Our goal in the present paper is to "generalize" simultaneously both of these results as follows:

\begin{thm}\label{new}
The operator 
$R_t:= \check T_t\circ T_{2t}\circ \check T_t$
has the following properties:
\begin{itemize}
\item
Regularization : $R_tf$ is $C^{1,1}$ for all bounded $f$ and all $t>0$.
\item
Approximation : If $f$ is uniformly continuous, then $R_tf$ converges uniformly to $f$ as $t\lto 0$.
\item
Pinching: If $u\geq v$ are two locally  bounded functions such that $u$ and $-v$ are $k$-semi-concave,
then the inequality 
$u\geq R_tf\geq v$ holds for each $t\leq k$ if $u\geq f\geq v$.
\end{itemize}
\end{thm}

Theorem \ref{new} does not, properly speaking, generalize Theorem \ref{LL}.
However, it offers a new (although similar) answer to the same problem:
approximating uniformly continuous functions on Hilbert spaces by $C^{1,1}$
functions with a simple explicit formula.

Because of its symmetric form, the regularizing operator $R_t$
enjoys some nicer properties than the Lasry-Lions operators.
For example, if $f$ is $C^{1,1}$, then it follows from the pinching property that 
$R_tf=f$ for $t$ small enough.

Theorem \ref{Il}, 
 can be  proved  using Theorem \ref{new} by taking $w=R_tu$, for $t$ small enough.
Note, in view of Lemma \ref{3} bellow, that $R_tu =\check T_t \circ T_t u$
when $t$ is small enough. 

Theorem \ref{new} can be somehow  extended to the case of 
finite dimensional open sets or
manifolds via partition of unity, at the price of loosing the simplicity of explicit expressions.
Let $M$ be a paracompact manifold of dimension $n$, equipped once and for 
all with an atlas $(\phi_i,i\in \Im)$ composed of charts
$\phi_i:B^n\lto M$, where $B^n$ is the open unit ball of radius 
one centered at the origin in $\Rm^n$.
We assume in addition that the image $\phi_i(B^n)$ is a relatively compact open set.
Let us fix, once and for all, a partition of the unity $g_i$ subordinated
to the open covering $(\phi_i(B^n),i\in \Im)$.
It means that the function $g_i$ is non-negative, with support inside $\phi_i(B^n)$,
such that $\sum_i g_i=1$, where the sum is locally finite.
Let us define the following formal operator 
$$
G_t(u):=
\sum_i 
\big[R_{ta_i}\big( (g_i u)\circ \phi_i
\big)
\big]\circ \phi_i^{-1},
$$
 where $a_i,i\in \Im$ are positive real numbers.
 We say that a function $u:M\lto \Rm$ is locally semi-concave if, for each 
 $i\in \Im$, there exists a constant $b_i$ such that the function 
 $u\circ \phi_i-\|.\|^2/{b_i}$ is concave on $B^n$.
 
 \begin{thm}\label{manifold}
Let  $u\geq v$ be two continuous functions on $M$ such that 
$u$ and $-v$ 
are locally semi-concave.
Then, the real numbers $a_i$ can be chosen such that,
 for each $t\in ]0,1]$ and each function $f$ satisfying $u\geq f\geq v$, we have:
\begin{itemize}
\item The sum in the definition of $G_t(f)$ is locally finite, so that the function
$G_t(f)$ is well-defined.
\item
The function $G_tf$ is locally $C^{1,1}$.
\item
If $f$ is continuous, then  $G_t(f)$ converges locally uniformly to $f$
as $t\lto 0$.
\item $u\geq G_tf\geq v$.
 \end{itemize}

 \end{thm}
 
 \subsection*{Notes and Acknowledgements}
 
Theorem \ref{Il} appears in Ilmanen's paper \cite{Il} as Lemma 4G. 
Several proofs are sketch there but none is detailed.
The proof we detail here follows lines similar to one of the sketches
of Ilmanen.
 This statement also has a more geometric counterpart, Lemma 4E
in \cite{Il}.
A detailed proof of this geometric version is given 
in \cite{C}, Appendix. 
My attention was attracted to these statements 
and their relations with recent progresses on sub-solutions 
of the Hamilton-Jacobi equation (see \cite {FS,ENS,Z})
by Pierre Cardialaguet, Albert Fathi
and Maxime Zavidovique.
 These authors also recently wrote a detailed proof
of Theorem \ref{Il}, see \cite{CFZ}.
This paper also proves  how  the geometric version 
follows from  Theorem \ref{Il}.
There are many similarities between the tools used in the present paper and
those used in \cite{ENS}.
Moreover, Maxime Zavidovique observed in \cite{Z} that the existence of $C^{1,1}$ subsolutions
of the Hamilton-Jacobi equation in the discrete case can be deduced from
Theorem \ref{Il}.
However, is seems that the main result of \cite{ENS}
(the existence of $C^{1,1}$ subsolutions in the continuous case) can't be deduced easily
from Theorem \ref{Il}. Neither can Theorem \ref{Il} be deduced from it.

\section{The operators $T_t$ and $\check T_t$ on Hilbert spaces}

The proofs  of the theorems  follow 
from standard properties of the operators $T_t$ and $\check T_t$
that we now recall in details.

\begin{lem}\label{1}
For each bounded function $u$, the function $T_tu$ is $t$-semi-concave and the function $\check T_tu$
is $t$-semi-convex.
Moreover, if $u$ is $k$-semi-concave, then for each $t<k$ the function 
$\check T_t u$ is $(k-t)$-semi-concave.
Similarly, if $u$ is $k$-semi-convex , then for each $t<k$ the function 
$T_t u$ is $(k-t)$-semi-convex.
\end{lem}
\proof
We shall prove the statements concerning $T_t$.
We have 
$$T_tu(x)-\|x\|^2/t=
\inf_y \big( u(y)+\|y-x\|^2/t-\|x\|^2/t\big)=
\inf_y \big( u(y)+\|y\|^2/t -2x\cdot y/t\big),
$$
this function is convex as an infimum of linear functions.
On the other hand, we have
$$T_tu(x)+\|x\|^2/l=
\inf_y \big( u(y)+\|y-x\|^2/t+\|x\|^2/l\big).
$$
Setting $f(x,y):= u(y)+\|y-x\|^2/t+\|x\|^2/l$,
the function $\inf_y f(x,y)$ is a convex function of $x$ if 
$f$ is a convex function of $(x,y)$.
This is true if $u$ is $k$-semi-convex, $t<k$, and  $l=k-t$ 
because we have the expression
$$f(x,y)=u(y)+\|y-x\|^2/t+\|x\|^2/l
=
(u(y)+\|y\|^2/k)+\big\|\sqrt{\frac{l}{kt}}y-\sqrt{\frac{k}{lt}}x
\big\|^2.
$$
\qed

Given a uniformly continuous function $u:H\lto \Rm$, we define
its modulus of continuity $\rho(r):[0,\infty)\lto [0,\infty)$
by the expression 
$\rho(r)=\sup_{x,e} u(x+re)-u(x),
$
where the supremum is taken on all $x\in H$ and all $e$ in the unit ball of $H$.
The function $\rho$ is non-decreasing, it satisfies $\rho(r+r')\leq \rho(r)+\rho(r')$,
and it converges to zero in zero (this last fact is equivalent to the uniform continuity of $u$).
We say that a function $\rho:[0,\infty)\lto [0,\infty)$ is a modulus of continuity
if it satisfies these properties.
Given a modulus of continuity $\rho(r)$, we say that a function $u$
is $\rho$-continuous if $|u(y)-u(x)|\leq \rho(\|y-x\|)$ for all $x$ and $y$ in  $H$.

\begin{lem}\label{2}
If $u$ is uniformly continuous,
then the functions $T_tu$ and $\check T_tu$ converge uniformly to $u$ when $t\lto 0$.
Moreover, given a modulus of continuity  $\rho$, there exists a non-decreasing function 
$\epsilon (t):[0,\infty)\lto [0,\infty)$ satisfying $\lim_{t\lto 0} \epsilon(t)=0$
and such that, for each $\rho$-continuous bounded function $u$, we have: 
\begin{itemize}
\item $T_tu$ and $\check T_t u$ are $\rho$-continuous for each $t\geq 0$.
\item $u-\epsilon(t)\leq T_tu(x)\leq u$ and 
$u\leq \check T_t u \leq u+\epsilon(t)$ for each $t\geq 0$.
\end{itemize}
\end{lem}

\proof
Let us fix $y\in H$, and set $v(x)=u(x+y)$.
We have $u(x)-\rho(|y\|)\leq v(x)\leq u(x)+\rho(\|y\|)$.
Applying the operator $T_t$ gives
$T_tu(x)-\rho(y)\leq T_tv(x)\leq 
T_tu(x)+\rho(y).
$
On the other hand, 
we have 
$$T_tv(x)=
\inf_z \big(u(z+y)+\|z-x)\|^2/t\big)=\inf _z\big( u(z)+\|z-(x+y)\|^2/t\big)
=T_tu(x+y),$$
so that 
$$T_tu(x)-\rho(\|y\|)\leq T_tu(x+y)\leq T_t u(x)+\rho(\|y\|).
$$
We have proved that $T_tu$ is $\rho$ continuous if $u$ is,
the proof for $\check T_tu$ is the same.

In order to study the convergence, let us set $\epsilon(t)=\sup_{r>0}(\rho(r)-r^2/t)$.
We have 
 $$
\epsilon(t) = \sup_{r> 0}\big( \rho(r\sqrt{t})-r^2\big)
\leq\sup_{r> 0}  \big((r+1)\rho(\sqrt{t})-r^2\big)
 =\rho(\sqrt{t})+\rho^2(\sqrt{t})/4.
$$
We conclude that $\lim_{t\lto 0} \epsilon(t)=0$.
We now come back to the operator $T_t$, 
and observe that
$$
u(y)-\|y-x\|^2/t \geq u(x)-\rho(\|y-x\|)+\|y-x\|^2/t\geq
u(x)-\epsilon(t)
$$ 
for each $x$ and $y$,
so that 
$$u-\epsilon(t)\leq T_tu \leq u.$$
\qed

\begin{lem}\label{3}
For each locally bounded function $u$, we
 have 
$\check T_t\circ T_t (u)\leq u
$ 
and the equality
$\check T_t \circ T_t (u)=u$
holds if and only if $u$ is   $t$-semi-convex.
Similarly, given a locally bounded function $v$,
we have $T_t\circ \check T_t (v)\geq v$, with equality if and only if
$v$ is $t$-semi-convex.
\end{lem}
\proof
 Let us write explicitly
 $$
 \check T_t \circ T_tu(x)=\sup_y\inf _z \big(u(z)+\|z-y\|^2/t-\|y-x\|^2/t\big).
 $$
 Taking $z=x$, we obtain the estimate 
 $
  \check T_t \circ T_tu(x)\leq \sup_y u(z)=u(z)
  $.
 Let us now write
  $$
 \check T_t \circ T_tu(x)+\|x\|^2/t
 =\sup_y\inf _z \big(u(z)+\|z\|^2/t+(2y/t)\cdot(x-z) \big)
$$
which by an obvious change of variable leads to
$$
\check T_t \circ T_tu(x)+\|x\|^2/t
 =\sup_y\inf _z \big(u(z)+\|z\|^2/t+y\cdot(x-z) \big).
$$
We recognize here  that the function $\check T_t \circ T_tu(x)+\|x\|^2/t$ is the
Legendre bidual of the function $u(x)+\|x\|^2/t$.
It is well-know that a locally bounded function is equal to its Legendre bidual
if and only if it is convex.
\qed

\begin{lem}\label{4}
If $u$ is locally bounded and semi-concave, then $\check T_t\circ T_t u$
is $C^{1,1}$ for each $t>0$.
\end{lem}

\proof
Let us assume that $u$ is $k$-semi-concave.
Then $u=T_k\circ \check T_ku$,  by Lemma \ref{3}.
We conclude that $\check T_t \circ T_t u = \check T_t \circ T_{t+k} f$,
with $f=\check T_k u$.
By Lemma \ref{1}, the function $T_{t+k}f$ is $(t+k)$-semi-concave.
Then, the function $\check T_t T_{t+k}f$ is $k$-semi-concave.
Since it is also $t$-semi-convex, it is $C^{1,1}$.
\qed

\section{Proof of the main results}
\textsc{Proof of Theorem \ref{new}:}
For each function $f$ and each $t>0$,
the function $\check T_t \circ T_{2t}\circ \check T_t f$
is $C^{1,1}$. This is a consequence of 
Lemma \ref{4} since 
$$\check T_t\circ T_{2t}\circ \check T_t f=\check T_t \circ T_t (T_t\circ \check T_t f)
$$
and since the function $T_t\circ \check T_t f$
is semi-concave.

Assume now that both $u$ and $-v$ are $k$-semi-concave.
We claim that $$
u\geq f\geq v \Longrightarrow u\geq  T_{t}\circ \check T_t f\geq v
\text{ and } u\geq \check T_t\circ T_{t} f\geq v
$$
for $t\leq 1/k$.
This claim implies that $u\geq \check T_t\circ T_{2t}\circ \check T_t f\geq v$
when $u\geq f\geq v$.
Let us now prove the claim concerning $\check T_t\circ T_t$, the other part being similar.
Since $v$ is $k$-semi-convex, we have $\check T_t\circ T_t v=v$
for $t\leq k$, by Lemma \ref{3}. Then,
$$u\geq f\geq \check T_t\circ T_tf \geq \check T_t\circ T_t v=v$$
where the second inequality follows from Lemma \ref{3}, and the third
from the obvious fact that the operators $T_t $ and $\check T_t$
are order-preserving.

The approximation property follows directly from Lemma \ref{2}.
\qed

\textsc{ Proof of Theorem \ref{manifold}:}
Let $a_i$ be chosen such that the functions 
 $(g_iu)\circ \phi_i$
and $-(g_iv)\circ \phi_i$ 
are $a_i$-semi-concave on $\Rm^n$.
The existence of real numbers $a_i$ with this property follows
from Lemma \ref{localization} below.
Given  $u\geq f\geq v$,
we can apply Theorem \ref{new} for each $i$ to the functions 
$$(g_iu)\circ \phi_i\geq (g_if)\circ \phi_i\geq (g_iv)\circ \phi_i
$$
 extended by zero outside of $B^n$.
We conclude that, for $t\in ]0,1]$, the function 
$R_{ta_i}((g_if)\circ \phi_i)$ is $C^{1,1}$
and satisfies 
$$ (g_iu)\circ \phi_i\geq R_{ta_i}((g_if)\circ \phi_i)
\geq (g_iv)\circ \phi_i.
$$
As a consequence, the function 
$$\big[R_{ta_i}\big((g_if)\circ \phi_i\big)\big]\circ \phi_i^{-1}
$$
is null outside of the support of $g_i$, and therefore the sum in the definition
of $G_t f$ is locally finite.
The function $G_t(f)$ is thus locally a finite sum
of $C^{1,1}$ functions hence it is locally $C^{1,1}$.
Moreover, we have 
$$
u=\sum_i g_iu\geq G_t(f)\geq \sum_ig_i v=v.
$$
\qed

We have used:
\begin{lem}\label{localization}
Let $u:B^n\lto \Rm$ be a bounded function
 such that $u-\|.\|^2/a$ is concave, for some $a>0$.
For each compactly supported non-negative  $C^2$ function $g:B^n\lto \Rm$,
the product $gu$ (extended by zero outside of $B^n$) is semi-concave on $\Rm^n$.
\end{lem}

\proof
Since $u$ is bounded, we can assume that  $u\geq 0$ on $B^n$.
Let $K\subset B^n$ be a compact subset of the open ball $B^n$
which contains the support of $g$ in its interior.
Since the function $u-\|.\|^2/a$ is concave on $B_1$
it admits super-differentials at each point.
As a consequence, for each $x\in B^n$, there exists a linear form $l_x$ such that 
$$0\leq u(y)\leq u(x)+l_x\cdot (y-x) +\|y-x\|^2/a
$$
for each $y\in B^1$.
Moreover, the linear form $l_x$ is bounded independently of $x\in K$.
We also have 
$$0\leq g(y)\leq g(x)+dg_x\cdot (y-x)+C\|y-x\|^2
$$
for some $C>0$, for all $x,y$ in $\Rm^n$.
Taking the product, we get, for $x\in K$ and $y\in B^n$,
$$u(y)g(y)\leq u(x)g(x)+(g(x)l_x+u(x)dg_x)\cdot (y-x)+
C\|y-x\|^2+C\|y-x\|^3 +C\|y-x\|^4
$$
where $C>0$ is a constant independent of $x\in K$ and $y\in B^n$, 
which may change from line to line.
As a consequence, setting $L_x=g(x)l_x+u(x)dg_x$, we obtain the inequality
\begin{equation}\tag{L}\label{L}
(gu)(y)\leq (gu)(x)+L_x\cdot (y-x)+C\|y-x\|^2
\end{equation}
for each $x\in K$ and $y\in B^n$.
If we set $L_x=0$ for $x\in \Rm^n-K$, the relation (\ref{L}) holds 
for each $x\in \Rm^n$ and $y\in \Rm^n$.
For $x\in K$ and $y\in B^n$, we have already proved it.
Since the linear forms $L_x$, $x\in K$ are uniformly bounded,
we can assume that $L_x\cdot (y-x)+C\|y-x\|^2\geq 0$ for 
all $x\in K$ and $y\in \Rm^n-B^n$ by taking $C$ large enough. Then, (\ref{L}) holds for all 
$x\in K$ and $y\in \Rm^n$.
For $x\in \Rm^n-K$ and $y$ outside of the support $g$,
 the relation (\ref{L}) holds in an obvious way, because
$gu(x)=gu(y)=0$, and $L_x=0$.
For $x\in \Rm^n-K$ and  $y$ in the support of $g$, the relation holds provided that 
$C\geq \max (gu)/d^2$, where $d$ is the distance between 
the complement of $K$ and the support of $g$. This is a positive number since $K$ 
is a compact set containing the support of $g$ in its interior.
We conclude that the function $(gu)$ is semi-concave on $\Rm^n$.
\qed

For completeness, we also prove, following Fathi:

\begin{lem}\label{F}
Let $u$ be a continuous function which is both $k$-semi-concave and $k$-semi-convex.
Then the function $u$ is $C^{1,1}$, and $6/k$ is a Lipschitz constant for the gradient of $u$.
\end{lem}
\proof
Let $u$ be a continuous function which is both $k$-semi-concave and $k$-semi-convex.
Then, for each $x\in H$, there exists a unique  $l_x\in H$ such that 
$$|u(x+y)-u(x)-l_x\cdot y|\leq \|y\|^2/k.$$
We conclude that $l_x$ is the gradient of $u$ at $x$, and we have to prove that the map
$x\lmto l_x$ is Lipschitz.
We have, for eah $x$, $y$ and $z$ in $H$:
$$
l_x\cdot (y+z)-\|y+z\|^2/k\leq u(x+y+z)-u(x)\leq l_x\cdot (y+z)+\|y+z\|^2/k
$$
$$
l_{(x+y)}\cdot (-y)-\|y\|^2/k\leq u(x)-u(x+y)\leq
l_{(x+y)}\cdot (-y)+\|y\|^2/k
$$
$$
l_{(x+y)}\cdot (-z)- \|z\|^2/k \leq u(x+y)-u(x+y+z)\leq l_{(x+y)}\cdot (-z) +\|z\|^2/k.
$$
Taking the sum, we obtain
$$
\big|(l_{x+y}-l_x)\cdot (y+z)
\big|\leq \|y+z\|^2/k+\|y\|^2/k+\|z\|^2/k.
$$
By a change of variables, we get 
$$
\big|(l_{x+y}-l_x)\cdot (z)
\big|\leq \|z\|^2/k+\|y\|^2/k+\|z-y\|^2/k.
$$
Taking $\|z\|=\|y\|$, we obtain
$$
\big|(l_{x+y}-l_x)\cdot (z)
\big|\leq 6\|z\|\|y\|/k
$$
for each $z$ such that $\|z\|=\|y\|$,
we conclude that 
$$\|l_{x+y}-l_x\|\leq 6\|y\|/k.
$$
\qed

\bibliographystyle{amsplain}
\providecommand{\bysame}{\leavevmode\hbox
to3em{\hrulefill}\thinspace}

\end{document}